\documentclass[a4paper, 10pt, conference]{ieeeconf}      

\IEEEoverridecommandlockouts                              

\usepackage{amsfonts}
\usepackage{amsmath}
\usepackage{hyperref,cite}
\usepackage{cleveref}
\usepackage{graphicx}
\usepackage{algorithm}
\usepackage{algorithmic}
\let\labelindent\relax
\usepackage{enumitem}
\usepackage{float}
\usepackage{caption}
\usepackage{booktabs}
\usepackage{makecell}
\usepackage{color}
\usepackage{threeparttable}
\usepackage{underscore}
\usepackage{multirow}

\newtheorem{theorem}{Theorem}
\DeclareMathOperator*{\argmin}{argmin} 

\newtheorem{lemma}[theorem]{Lemma}

\newcommand{\diag}{\mathop{\rm diag}\nolimits}
\renewcommand{\mod}{\mathop{\rm mod}\nolimits}
\newcommand{\rr}{{\mathbb R}}
\newcommand{\smallmat}[1]{\left[ \begin{smallmatrix}#1 \end{smallmatrix} \right]}

\begin{document}

\title{\LARGE \bf
A Simple and Fast Coordinate-Descent Augmented-Lagrangian Solver for Model Predictive Control}
\author{Liang Wu$^{1}$, Alberto Bemporad$^{1}$%
\thanks{The authors are with the IMT School for Advanced Studies Lucca, Italy,
        {\tt\small \{liang.wu,alberto.bemporad\}@imtlucca.it}}%
}
\thispagestyle{empty}
\pagestyle{empty}
\maketitle

\begin{abstract}
This paper proposes a novel Coordinate-Descent Augmented-Lagrangian (CDAL) solver for linear, possibly parameter-varying, model predictive control (MPC) problems. At each iteration, an augmented Lagrangian (AL) subproblem is solved by coordinate descent (CD), exploiting the structure of the MPC problem. The CDAL solver enjoys three main properties: ($i$) it is \emph{construction-free}, in that
it avoids explicitly constructing the quadratic programming (QP) problem associated with MPC;
($ii$) is \emph{matrix-free}, as it avoids multiplications and factorizations of matrices;
and ($iii$) is \emph{library-free}, as it can be simply coded without any library dependency, 
90-lines of C-code in our implementation. To favor convergence speed, CDAL employs a reverse cyclic rule for the CD method, the accelerated Nesterov’s scheme for updating the dual variables, a simple diagonal preconditioner, and an efficient coupling scheme between the CD and AL methods. We show that CDAL competes with other state-of-the-art methods, both in case of unstable linear time-invariant
and linear parameter-varying prediction models.
\end{abstract}


\section{Introduction}
Model predictive control (MPC) has been widely used for decades to control multivariable systems subject to input and output constraints~\cite{qin2003survey}. Apart from small-scale linear time-invariant (LTI) MPC problems whose explicit MPC control law can be obtained \cite{bemporad2002explicit}, deploying an MPC controller in an electronic control unit requires an embedded Quadratic Programming (QP) solver. In the past decades, the MPC community has made tremendous research efforts to develop embedded QP algorithms~\cite{KFZD18}, based on interior-point methods~\cite{wang2009fast,Wri19}, active-set algorithms~\cite{ferreau2008online,Bem16}, gradient projection methods~\cite{patrinos2013accelerated}, the alternating direction method of multipliers (ADMM)~\cite{boyd2011distributed,SBGBB20}, and other techniques~\cite{LS97,HTP19,Bem18,SB20,saraf2019efficient}.

A demanding requirement for industrial MPC applications is code simplicity, for easily being verified, validated, and maintained on embedded platforms. In this respect, the interior-point and active-set methods require more complicated arithmetic operations in their algorithm implementations when compared to first-order optimization methods like gradient projection and ADMM. The first-order optimization methods are quite appealing in embedded MPC since their embedded implementations could only involve additions and multiplications (no divisions, square roots, etc.). However, most of the proposed approaches require that the MPC-to-QP transformation is explicitly constructed for consumption by the solver, such as for preconditioning, estimating the Lipschitz constant of the cost gradient, and factorizing matrices. This may not be an issue for linear time-invariant (LTI) MPC problems, in which the MPC-to-QP construction and other operations on the problem matrices can be done off-line. But for some linear parameter-varying (LPV) or for linear time-varying MPC problems in which the linear dynamic model, cost function and/or constraints change at run time, an explicit online MPC-to-QP construction increases the complexity of the embedded code and computation time. Avoiding an explicit MPC-to-QP construction, can be called as construction-free property of an MPC solver. The barrier interior-point FastMPC solver~\cite{wang2009fast} and the active-set based BVLS solver~\cite{saraf2019efficient} are construction-free; they directly use the model and weight matrices to define the MPC problem without constructing a QP problem. Their complicated implementations are not matrix-free as involving Cholesky or QR factorizations arithmetic operations during iterations. The well-known simple and efficient first-order method OSQP~\cite{SBGBB20} is not construction-free and matrix-free when applied to solve LPV-MPC problems, as it requires that matrix factorizations are computed and cached on each sampling time. The OSQP utilizes its own $LDL^T$ solver to perform matrix factorizations, thus being library-free.

\subsection{Contribution}
By combining the coordinate descent (CD) and augmented Lagrangian (AL) methods, in this paper we develop a construction-free, matrix-free, and library-free solver for LTI and LPV MPC problems that is particularly suitable for embedded industrial deployment.

Coordinate descent has received extensive attention in recent years due to its application to machine learning \cite{hsieh2008dual, chang2008coordinate, richtarik2016distributed}. In this paper, we will exploit the special structure arising from linear MPC formulations when applying CD. In \cite{richter2011computational,  nedelcu2014computational, kogel2011fast}, the authors also use AL to solve linear MPC problems with input and state constraints using the fast gradient method \cite{Nes83} to solve the associated subproblems. The Lipschitz constant of the cost gradient and convexity parameters \cite{richter2011computational} are needed to achieve convergence, and computing them requires in turn the Hessian matrix of the subproblem, and hence constructing the QP problem. As the Hessian matrix of the AL subproblem is close to a block diagonal matrix, this suggests the use of the CD method to solve such a QP subproblem, due to the fact that CD does not require any problem-related parameter. Moreover, only small matrices are involved in running the CD method, namely the matrices of the linear prediction model and the weight matrices. As a result, the proposed CDAL algorithm does not require the QP construction phase and is extremely simple to implement. In addition, each update of the optimization vector has a computation cost per iteration that is quadratic with the state and input dimensions and linear with the prediction horizon.

To improve the convergence speed of CDAL, we propose four techniques: a reverse cyclic rule
for CD, Nesterov's acceleration~\cite{Nes83}, preconditioning, and an efficient coupling between CD and AL. While the use of a reverse cyclic rule in CD still preserves  convergence, when the MPC problem is solved by warm-starting it from the shifted previous optimal solution, the gap between the initial guess and the new optimal solution is mainly caused by the last block of variables, and computing the last block at the beginning tends to reduce the overall number of required iterations to converge, as we will verify in the numerical experiments reported in this paper. We employ Nesterov's acceleration scheme for updating the dual vector to improve computation speed and a heuristic preconditioner that simply scales the state variables. In addition, an efficient coupling scheme between CD and AL method is proposed to reduce the computation cost of each CD iteration. 
To analyze the role of each component of CDAL and its computational performance with respect to other solvers (FastMPC, $\mu$AO-MPC, OSQP, and MATLAB's \texttt{quadprog}), we conduct numerical experiments on an ill-conditioned problem of LTI-MPC control of an open-loop unstable AFTI-16 aircraft, and on LPV-MPC control of a continuously stirred tank reactor (CSTR).


\subsection{Notation}
$H \succ 0$ ($H \succeq 0$) denotes positive definiteness (semi-definiteness) of a square matrix $H$, $H^{\prime}$ (or $z^{\prime}$) denotes the transpose of matrix $H$ (or vector $z$), $H_{i,j}$ denotes the element of matrix $H$ on the $i$th row and the $j$th column, $H_{i,\cdot}, H_{\cdot,j}$ denote the $i$th row vector, and  $j$th column vector of matrix $H$, respectively. For a vector $z$, $\|z\|_2$ denotes the Euclidean norm of $z$, $z_{\neq i}$ the subvector obtained from $z$ by eliminating its $i$th component $z_i$. 

\section{Model Predictive Control}
\label{sec:MPC}
Consider the following MPC formulation for tracking problems
\begin{eqnarray}
\min && \frac{1}{2}\sum_{t=0}^{T-1}\left\|W_y \left(y_{t+1}-r_{t+1}\right)\right\|_{2}^{2} + \frac{1}{2}\left\|W_u  \left(u_{t}-u^r_{t}\right)\right\|_{2}^{2} \nonumber\\&&\hspace*{2cm}+ \frac{1}{2}\left\|W_{\Delta u}  \Delta u_{t}\right\|_{2}^{2} \nonumber\\
\text {s.t.}  && x_{t+1}=A x_{t}+B u_{t},\ t=0, \ldots, T-1\nonumber\\
&& y_{t+1} = C x_{t+1},\ t=0, \ldots, T-1\nonumber\\
&& u_{t} = u_{t-1} + \Delta u_{t},\ t=0, \ldots, T-1\nonumber\\
&& x_{\min } \leq x_{t} \leq x_{\max},\ t=1, \ldots, T\nonumber\\
&& u_{\min } \leq u_{t} \leq u_{\max},\ t=0, \ldots, T-1\nonumber\\
&& \Delta u_{\min } \leq \Delta u_{t} \leq \Delta u_{\max }, t= 0, \ldots, T-1\nonumber\\
&& x_0 = \bar{x}_0, u_{-1} = \bar{u}_{-1}
\label{MPC}
\end{eqnarray}
in which $x_t\in\rr^{n_x}$ is the state vector, $u_t\in\rr^{n_u}$ the input vector,
$\Delta u_{t}=u_{t}-u_{t-1}$ the vector of input increments,
$y_t\in\rr^{n_y}$ the output vector, $r_t$ and $u_t^r$ are the output and input set-points, 
and $\bar{x}_0$ and $\bar{u}_{-1}$ denote the current state and the previous input vectors, respectively. We assume that $W_y=W_y^\prime\succeq 0$,
$W_u=W_u^\prime\succeq 0$, $W_{\Delta u}=W_{\Delta u}^\prime\succ 0$. The formulation~\eqref{MPC} could
be extended to include time-varying bounds on $x$ and $u$ along the prediction horizon, 
linear equality constraints or box constraints on the terminal state $x_T$
for guaranteed closed-loop convergence, as well as affine prediction models. To simplify the notation, in the sequel
we consider the following reformulation of~\eqref{MPC}
\begin{eqnarray}
\min && \frac{1}{2}\sum_{t=1}^{T}\hat{x}_{t}^{\prime}(\hat{C}^{\prime}\hat{W}\hat{C})\hat{x}_{t} - \hat{x}_{t}^{\prime}(\hat{C}^{\prime}\hat{W}\hat{r}_{t})
+ \frac{1}{2}\hat{u}_{t-1}^{\prime}W_{\Delta u} \hat{u}_{t-1}\nonumber\\
\text {s.t.}  && \hat{x}_{t+1}=\hat{A} x_{t}+\hat{B} \hat{u}_{t}, \ t=0, \ldots, T-1\nonumber\nonumber\\
&& \hat{x}_{\min} \leq \hat{x}_{t} \leq \hat{x}_{\max},\ t=1, \ldots, T\nonumber\\
&& \hat{u}_{\min} \leq \hat{u}_{t} \leq \hat{u}_{\max},\ t=0, \ldots, T-1\nonumber\\
&& \hat{x}_0 = \smallmat{\bar{x}_0\\ \bar{u}_{-1}}\label{MPC2}
\end{eqnarray}
where $\hat{x}_{t} = \smallmat{x_t\\u_{t-1}} \in \mathcal{R}^{\hat{n}_x}$, 
$\hat{n}_x=n_x+n_u$, $\hat{u}_t = \Delta u_t \in \mathcal{R}^{n_u}$, $\hat{A}=\smallmat{
A & B \\
0 & I
}\in \mathbb{R}^{\hat{n}_x \times \hat{n}_x}$, $\hat{B} = \smallmat{
B \\
I
}\in \mathbb{R}^{\hat{n}_x \times n_u}$, $\hat{C}=\smallmat{
C & 0 \\
0 & I}$, $\hat{W} = \smallmat{
W_y & 0 \\
0 & W_u}$, $\hat{r}_t = \smallmat{
r_t \\
u^r_{t-1}
}$.
The vector $z$ of variables to optimize is
\[
z=\left[\begin{array}{cccccc}
\hat{u}_{0}^{\prime} & \hat{x}_{1}^{\prime} & \hat{u}_{1}^{\prime} & \ldots & \hat{u}_{T-1}^{\prime} & \hat{x}_{T}^{\prime}
\end{array}\right]^{\prime}\in\rr^{T(\hat{n}_x+n_u)}
\]
The inequality constraints on state and input variables, whose number is $2T(\hat{n}_x+n_u)$, are 
\[
\underline{z} \leq z \leq \bar{z} \Leftrightarrow\left\{\begin{array}{l}
\hat{x}_{\min} \leq \hat{x}_{t} \leq \hat{x}_{\max}, \forall t=1, \ldots, T \\
\hat{u}_{\min} \leq \hat{u}_{t} \leq \hat{u}_{\max}, \forall t=0, \ldots, T-1
\end{array}\right.
\]
where $\hat{x}_{\min} =  \smallmat{x_{\min}\\u_{\min}}$, $\hat{x}_{\max} =  \smallmat{x_{\max}\\u_{\max}}$, $\hat{u}_{\min} =  \Delta u_{\min}$ and $\hat{u}_{\max} =  \Delta u_{\max}$. At each sample step, the MPC problem~\eqref{MPC} can be recast as the following quadratic program (QP)
\begin{eqnarray}
\min && \frac{1}{2} z^{\prime} H z+h^{\prime} z \nonumber\\
\text { s.t. } && \underline{z} \leq z \leq \bar{z} \nonumber\\
&& G z=g\label{QP}
\end{eqnarray}
where $H =H^\prime\succeq 0$, $H\in \mathbb{R}^{n_z \times n_z}$, $n_z=T(\hat{n}_x+n_u)$, $h \in \mathbb{R}^{n_z}$, $G \in \mathbb{R}^{T\hat{n}_x \times n_z}$, and $g \in \mathbb{R}^{T\hat{n}_x}$ are defined as
\begin{eqnarray*}
H&=&{\small \left[\begin{array}{ccccc}
R & 0 & \ldots & 0 & 0 \\
0 & Q & \ldots & 0 & 0 \\
\vdots & \vdots & \ddots & \vdots & \vdots \\
0 & 0 & \ldots & R & 0 \\
0 & 0 & \ldots & 0 & Q
\end{array}\right]},\ \begin{array}{rcl}
    R &=& W_{\Delta u}\\~\\
    Q &=& \hat{C}^{\prime}\hat{W}\hat{C}\end{array}\\
G&=&{\small \left[\begin{array}{cccccccc}
\hat{B} & -I & 0 & 0 & \ldots & 0 & 0 & 0 \\
0 & \hat{A} & \hat{B} & -I & \ldots & 0 & 0 & 0 \\
\vdots & \vdots & \vdots & \vdots & \ddots & \vdots & \vdots & \vdots \\
0 & 0 & 0 & 0 & \ldots & \hat{A} & \hat{B} & -I
\end{array}\right]}\\
h&=&{\small \left[\begin{array}{c}
-\hat{C}^{\prime}\hat{W}\hat{r}_{1} \\
-\hat{C}^{\prime}\hat{W}\hat{r}_{2} \\
\vdots \\
-\hat{C}^{\prime}\hat{W}\hat{r}_{T} \\
\end{array}\right]},\ 
g={\small \left[\begin{array}{c}
-\hat{A}\hat{x}_{0}\\
0 \\
\vdots \\
0
\end{array}\right]}
\end{eqnarray*}
Clearly matrix $G$ is full row-rank.
Note that $A,B,C,W_y,W_u,W_{\Delta u}$ and the upper and lower bounds on $x$, $u$, 
and $\Delta u$ in~\eqref{MPC} may change at each controller execution. 

\section{Algorithm}
\label{sec:algo}
\subsection{Augmented Lagrangian Method}
We solve the convex quadratic programming problem~\eqref{QP}
by applying the augmented Lagrangian method. The bound-constrained Lagrangian function $\mathcal{L}: \mathcal{Z} \times \mathbb{R}^{T \times \hat{n}_x} \rightarrow \mathbb{R}$ is given by
\[
\mathcal{L}(z,\Lambda) = \frac{1}{2} z^{\prime} H z + z^{\prime}h  + \Lambda^{\prime}(Gz-g)
\]
where $\mathcal{Z}= \left\{\underline{z} \leq z \leq \bar{z}\right\}$ and $\Lambda \in \mathbb{R}^{T\hat{n}_x}$ is the vector of Lagrange multipliers associated with the equality constraints in~\eqref{QP}. 
The dual problem of~\eqref{QP} is
\begin{equation}\label{DualQP}
\max_{\Lambda \in \mathbb{R}^{T\hat{n}_x}} \phi(\Lambda)
\end{equation}
where $\phi(\Lambda) = \min_{z\in \mathcal{Z}} \mathcal{L}(z,\Lambda)$. Assuming that Slater’s constraint qualification holds, the optimal value of the primal problem~\eqref{QP} and of its dual~\eqref{DualQP}
coincide. However, $\phi(\Lambda)$ is not differentiable in general~\cite{bertsekas1997nonlinear},
so that any subgradient method for solving~\eqref{DualQP} would have a slow convergence rate. 
Under the AL framework, the augmented Lagrangian function
\begin{equation}\label{AugLagFun}
\mathcal{L}_{\rho}(z,\Lambda) = \frac{1}{2} z^{\prime} H z+z^{\prime}h  + \Lambda^{\prime}(Gz-g) + \frac{\rho}{2} \|Gz-g\|^{2}
\end{equation}
is used instead, where the parameter $\rho > 0$ is a penalty parameter. The corresponding augmented dual problem is defined as:
\begin{equation}\label{AugDualQP}
\max_{\Lambda \in \mathbb{R}^{T \times n_x}} \phi_{\rho}(\Lambda)
\end{equation}
where $\phi_{\rho}(\Lambda) = \min_{z\in \mathcal{Z}} \mathcal{L}_{\rho}(z,\Lambda)$ is differentiable provided that $H+\rho G^{\prime}G \succ 0$. The dual problem ~\eqref{DualQP} and the augmented dual problem~\eqref{AugDualQP} share the same optimal solution~\cite[see chapter 2 subsection 2.2]{bertsekas2014constrained}, and most important $d_{\rho}(\Lambda)$ is concave and differentiable, with gradient~\cite{bertsekas1997nonlinear, nesterov2005smooth}
$\nabla \phi_{\rho}(\Lambda)=G z^{*}(\Lambda)-g$,
where $z^{*}(\Lambda)$ denotes the optimal solution of the inner problem $\min_{z\in \mathcal{Z}} \mathcal{L}_{\rho}(z,\Lambda)$ for a given $\Lambda$.
Moreover, the gradient mapping $\nabla \phi_{\rho}:\rr^{T \times n_x}\to\rr^{T \times n_x}$ is Lipschitz continuous, with a Lipschitz constant given by
$L_{\phi}=\rho^{-1}$~\cite{lan2016iteration}. 

Let $F_{\rho}(z;\Lambda^k) = \frac{1}{2} z^{\prime} H_A z+(h_A^k)^{\prime} z$,
where $h_A^k=\frac{1}{\rho}h+G^{\prime}\Lambda^k-G^{\prime}g$,
and $H_A=\frac{1}{\rho}H+G'G$ has the block-sparse structure
\[
H_{A}={\small\left[\begin{array}{ccccccccc}
\phi_{1} & \phi_{2} & 0 & 0 & 0 & \ldots & 0 & 0 & 0 \\
\phi_{2}^{\prime} & \phi_{3} & \phi_{4} & \phi_{5} & 0 & \ldots & 0 & 0 & 0 \\
0 & \phi_{4}^{\prime} & \phi_{1} & \phi_{2} & 0 & \ldots & 0 & 0 & 0 \\
0 & \phi_{5}^{\prime} & \phi_{2}^{\prime} & \phi_{3} & \phi_{4} & \ldots & 0 & 0 & 0 \\
\vdots & \vdots & \vdots & \vdots & \vdots & \ddots & \vdots & \vdots & \vdots \\
0 & 0 & 0 & 0 & 0 & \ldots & \phi_{3} & \phi_{4} & \phi_{5} \\
0 & 0 & 0 & 0 & 0 & \ldots & \phi_{4}^{\prime} & \phi_{1} & \phi_{2} \\
0 & 0 & 0 & 0 & 0 & \ldots & \phi_{5}^{\prime} & \phi_{2}^{\prime} & \phi_{6}
\end{array}\right]}
\]
and $\phi_{1}=\frac{1}{\rho}R+ \hat{B}^{\prime}\hat{B}$, $\phi_{2}=-\hat{B}^{\prime}$,
$\phi_{3}=\frac{1}{\rho}Q+\left(I+\hat{A}^{\prime}\hat{A}\right)$, $\phi_{4}=\hat{A}^{\prime}\hat{B}$,
$\phi_{5}=-\hat{A}^{\prime}$, $\phi_{6}=\frac{1}{\rho}Q+ I$. 
Since $G$ is full rank, matrix $H_A \succ 0$. According to \cite{bertsekas2014constrained}, the AL algorithm can be formulated in scaled form as follows:
\begin{subequations}\label{ALM}
\begin{eqnarray}
z^{k+1} &=& \argmin_{z \in \mathcal{Z}}  F_{\rho}(z;\Lambda^k)\label{subproblem2}\label{subproblem}\\
\Lambda^{k+1} &=& \Lambda^{k} + (Gz^{k+1}-g)
\end{eqnarray}
\end{subequations}
which involves the minimization step of the primal vector $z$ and the update step of the dual vector
$\Lambda$. 
As shown in~\cite{bertsekas2014constrained}, the convergence of AL can be assured for a large range of values of $\rho$. Typically, the larger the penalty parameter, the faster the AL algorithm is
to converge, but the more difficult~\eqref{subproblem} is to solve, due to a larger condition number of the Hessian matrix of subproblem~\eqref{subproblem}. The convergence rate of the AL algorithm~\eqref{ALM} is $O(1/k)$ according to \cite{he2010acceleration}.
To improve the speed of the AL method,  \cite{kang2015inexact} proposed an accelerated AL algorithm, whose iteration-complexity is $O(1/k^2)$ for linearly constrained convex programs, by using Nesterov's acceleration technique. The accelerated AL algorithm is summarized in Algorithm~\ref{AALM}.
\begin{algorithm}[H]
    \caption{Accelerated augmented Lagrangian method \cite{kang2015inexact}}\label{AALM}
    \textbf{Input}:  Initial guess $z^0 \in \mathcal{Z}$ and $\Lambda^0$;
maximum number $N_{\rm out}$ of iterations; parameter $\rho>0$.
    \vspace*{.1cm}\hrule\vspace*{.1cm}
    \begin{enumerate}[label*=\arabic*., ref=\theenumi{}]
    \item Set $\alpha_1\leftarrow1$; $\hat{\Lambda}^0\leftarrow\Lambda^0$;
        \item \textbf{for} $k=1, 2,\cdots, N_{\rm out}$ \textbf{do}
        \begin{enumerate}[label=\theenumi{}.\arabic*., ref=\theenumi{}.\arabic*]
            \item $z^{k} \leftarrow \argmin_{z \in \mathcal{Z}}  F_\rho(z;\hat{\Lambda}^{k-1})$;       
            \item $\Lambda^{k} \leftarrow \hat{\Lambda}^{k-1} + (G z^{k}-g)$;   
            \item \textbf{if} $\|\Lambda^k-\hat{\Lambda}^{k-1} \|_2^2 \leq \epsilon$, \textbf{stop};
            \item $\alpha_{k+1} \leftarrow \frac{1+\sqrt{1+4\alpha_{k}^2}}{2} $;
            \item $\hat{\Lambda}^{k} \leftarrow \Lambda^k + \frac{\alpha_{k}-1}{\alpha_{k+1}} (\Lambda^k-\Lambda^{k-1})$;
        \end{enumerate}
    \item \textbf{end}.
    \end{enumerate}
\end{algorithm}

For solving the strongly convex box-constrained QP~\eqref{subproblem2}, the fast gradient projection method was used in \cite{richter2011computational,kogel2011fast}. Inspired by the fact that the Gauss-Seidel method in solving block tridiagonal linear systems is efficient~\cite{amodio1995parallel}, in this paper we propose the use of the cyclic CD method to make full use of block sparsity and avoid the explicit construction of matrix $H_A$. Note that in the gradient projection method or fast gradient projection method \cite{kogel2011fast}, the Lipschitz constant parameter deriving from matrix $H_A$ needs to be calculated or estimated to ensure convergence. Therefore, for linear MPC problems that change at runtime such methods would be less preferable
than cyclic CD. 
In this paper, by making full use of the structure of the subproblem, we will implement a cyclic CD 
method that requires less computations, as we will detail in the next section. 

\subsection{Coordinate Descent Method}
The idea of the CD method is to minimize the objective function along only one coordinate direction at each iteration, while keeping the other coordinates fixed~\cite{wright2015coordinate}. In \cite{luo1992convergence}, the authors showed that the CD method is convergent in convex differentiable minimization problems, and the rate of convergence is at least linear. We first give a brief introduction of the CD method to solve~\eqref{subproblem2}. Under the assumption that the set of optimal solutions is nonempty and that the objective function $F_\rho$ is convex, continuously differentiable, and strictly convex with respect to each coordinate, the CD method proceeds iteratively for $k=0,1,\ldots,$ as follows:
\begin{subequations}\label{CDM}
\begin{eqnarray}
\text{choose}~i_k \in \left\{1,2,\ldots,n_z\right\} \\
z_{i_k}^{k+1} = \argmin_{z_{i_k} \in \mathcal{Z}} F_\rho(z_{i_k},z_{\neq i_k}^k;\hat{\Lambda}^k)
\end{eqnarray}
\end{subequations}
where with a slight abuse of notation we denote by $F_\rho(z_{i_k},z_{\neq i_k}^k;\hat{\Lambda}^k)$
the value $F_\rho(z;\hat{\Lambda}^k)$ when $z_{\neq i_k}=z_{\neq i_k}^k$ is fixed.
The convergence of the iterations in (\ref{CDM}) for $k\rightarrow\infty$ depends on the rule used to choose the coordinate index $i_k$. In \cite{luo1992convergence}, the authors show that the \emph{almost cyclic rule} and \emph{Gauss-Southwell rule} guarantee convergence. Here we use the almost cyclic rule, that provides convergence according to the following lemma:

\begin{lemma}[\cite{luo1992convergence}]
\label{lemma:cd}
Let $\left\{z^k\right\}$ be the sequence of coordinate-descent iterates~\eqref{CDM}, where every coordinate index is iterated upon at least once on every $N$ successive iterations, $N\geq n_z$. The sequence $\left\{z^k\right\}$ converges at least linearly to the optimal solution $z^{*}$
of problem~\eqref{subproblem2}.
\end{lemma}
In this paper we will use the \emph{reverse cyclic rule}
\[
    i_k=n_z-(k\mod n_z)
\]
to exploit the fact that the shifted previous optimal solution
is used as a warm start. The chosen rule clearly satisfies the assumptions 
of Lemma~\ref{lemma:cd} for convergence.
The implementation of one pass through all $n_z$ coordinates using reverse cyclic CD is reported in Procedure~\ref{RCCD}. In Procedure~\ref{RCCD}, the Lagrangian variable $\hat{\Lambda} \in \mathbb{R}^{T \times \hat{n}_x}$ is divided into $\{\hat{\lambda}_0,\ldots,\hat{\lambda}_{t-1},\ldots,\hat{\lambda}_{T-1}\}$, where $\hat{\lambda}_{t-1} \in \mathbb{R}^{\hat{n}_x}$. For a given symmetric $M\in\rr^{n_s\times n_s} \succeq 0$, $d\in\rr^{n_s}$, the operator $\operatorname{CCD}_{[\underline{s},\bar{s}]} \left\{M,d\right\}$ used in Procedure \ref{RCCD} represents 
one pass iteration of the reverse cyclic CD method through
all $n_s$ coordinates $s_{n_s},\ldots,s_1$ for the following box-constrained QP
\begin{equation}
\min _{s \in [\underline{s},\bar{s}]} \frac{1}{2} s^{\prime}Ms +s^{\prime} d
\end{equation}
that is to execute the following $n_s$ iterations
\begin{equation}\label{BoxCDM}
\begin{array}{l}
\text { for } i=n_s,\ldots,1 \\
\quad \quad s_i \leftarrow \left[s_i-\frac{1}{M_{i,i}}(M_{i,\cdot}s+d_i)\right]_{\underline{s}_i}^{\bar{s}_i} \\
\text { end }
\end{array}
\end{equation}
where $\left[s_i\right]_{\underline{s}_i}^{\bar{s}_i}$ is the projection operator
\begin{equation}
\left[s_i\right]_{\underline{s}_i}^{\bar{s}_i} = \left\{\begin{array}{lll}
\bar{s}_i & \mbox{if} & s_i \geq \bar{s}_i\\
s_i & \mbox{if} &\underline{s}_i < s_i < \bar{s}_i\\
\underline{s}_i&\mbox{if} & s_i \leq \underline{s}_i
\end{array}\right.
\end{equation}
Note that in Procedure~\ref{RCCD}, Steps~\ref{algo:rcd:ccd-1}, \ref{RCCD-stepT}, \ref{algo:rcd:ccd-2}, and \ref{RCCD-step} all involve the same operator CCD. In Procedure \ref{DetailedRCCD}, we exemplify an efficient way to evaluate such an operator for Step~\ref{RCCD-step} of Procedure \ref{RCCD}, as the approach is similar for evaluating Steps~\ref{algo:rcd:ccd-1}, \ref{RCCD-stepT}, and~\ref{algo:rcd:ccd-2}, where $\sigma$ records the sum of squared coordinate variations.

\floatname{algorithm}{Procedure}
\begin{algorithm}[t]
    \caption{Full pass of reverse cyclic coordinate descent on all block variables}\label{RCCD}
    \textbf{Input}: $\hat{\Lambda}=\{\hat{\lambda}_0,\ldots,\hat{\lambda}_{T-1}\}$, $U = \{\hat{u}_{0},\cdots, \hat{u}_{T-1}\}$, $X =\{\hat{x}_{0}, \hat{x}_{1}, \cdots, \hat{x}_{T}\}$; MPC settings $\hat{A}$, $\hat{B}$, $Q$, $R$, $\hat u_{\min}$, $\hat u_{\max}$, $\hat x_{\min}$, $\hat x_{\max}$; parameter $\rho>0$.
    \vspace*{.1cm}\hrule\vspace*{.1cm}
    \begin{enumerate}[label*=\arabic*., ref=\theenumi{}]
        \item $\sigma\leftarrow 0$;
        \item $\{\hat{x}_{T},\sigma\} \leftarrow \hspace*{-.3cm}\underset{\hat{x}_{T} \in [\hat x_{\min},\hat x_{\max}]}{\operatorname{CCD}} \{\frac{1}{\rho}Q+I, -\hat{\lambda}_{T-1}-\hat{A} \hat{x}_{T-1}-\hat{B} \hat{u}_{T-1}-\hat{C}^{\prime}\hat{W}\hat{r}_{T} , \sigma\}$;\label{algo:rcd:ccd-1}
        \item $\{\hat{u}_{T-1},\sigma\} \leftarrow \hspace*{-.3cm}\underset{\hat{u}_{T-1} \in [\hat u_{\min},\hat u_{\max}]}{\operatorname{CCD}} \{\frac{1}{\rho}R+\hat{B}^{\prime}\hat{B}, \hat{B}^{\prime}(\hat{\lambda}_{T-1}+\hat{A} \hat{x}_{T-1}-\hat{x}_{T}), \sigma \}$;
        \label{RCCD-stepT}
        \item \textbf{for} $t=T-2,T-3,\ldots,0$ \textbf{do}
        \begin{enumerate}[label=\theenumi{}.\arabic*., ref=\theenumi{}.\arabic*]
            \item $\{\hat{x}_{t+1}, \sigma\} \leftarrow \hspace*{-.3cm}\underset{\hat{x}_{t+1} \in [\hat x_{\min},\hat x_{\max}]}{\operatorname{CCD}} \{\frac{1}{\rho}Q+I+\hat{A}^{\prime} \hat{A}, -(\hat{\lambda}_{t}+\hat{A} \hat{x}_{t}+\hat{B} \hat{u}_{t}) + \hat{A}^{\prime}(\hat{\lambda}_{t+1}+\hat{B} \hat{u}_{t+1}-\hat{x}_{t+2} )-\hat{C}^{\prime}\hat{W}\hat{r}_{t} , \sigma \}$;
            \label{algo:rcd:ccd-2}
            \item $\{\hat{u}_{t}, \sigma\} \leftarrow \hspace*{-.3cm}\underset{\hat{u}_{t} \in [\hat u_{\min},\hat u_{\max}]}{\operatorname{CCD}} \{\frac{1}{\rho} R+\hat{B}^{\prime} \hat{B}, \hat{B}^{\prime}(\hat{\lambda}_{t}+\hat{A} \hat{x}_{t}-\hat{x}_{t+1}), \sigma \}$;\label{RCCD-step}
        \end{enumerate}
    \item \textbf{end}.
    \end{enumerate}
    \vspace*{.1cm}\hrule\vspace*{.1cm}
    \textbf{Output}: $\hat U$, $\hat X$, $\sigma$.
\end{algorithm}

\floatname{algorithm}{Procedure}
\begin{algorithm}[t]
    \caption{Evaluation of $\operatorname{CCD}$ in Step~\ref{RCCD-step} of Procedure \ref{RCCD}}\label{DetailedRCCD}
    \textbf{Input}: $\hat{\lambda}_t$, $\hat{u}_{t}$, $\hat{x}_{t}$, $\hat{x}_{t+1}$; MPC settings $\hat{A}$, $\hat{B}$, $R$, $\hat u_{\min}$, $\hat u_{\max}$; parameter $\rho>0$; update amount $\sigma\geq 0$.
    \vspace*{.1cm}\hrule\vspace*{.1cm}
    \begin{enumerate}[label*=\arabic*., ref=\theenumi{}]
        \item $V \leftarrow \hat{\lambda}_t + \hat{A}\hat{x}_{t} + \hat{B} \hat{u}_{t} - \hat{x}_{t+1}$;
        \item \textbf{for} $i=n_u,\ldots, 1$ \textbf{do}
        \begin{enumerate}[label=\theenumi{}.\arabic*., ref=\theenumi{}.\arabic*]
            \item $s \leftarrow  \frac{1}{\rho}R_{i,\cdot} \hat{u}_{t} + (\hat{B}_{\cdot, i})^{\prime} V$;
            \item $\theta \leftarrow  \left[\hat{u}_{t,i}-\frac{s}{\frac{1}{\rho}R_{ii}+(\hat{B}^{\prime}\hat{B})_{ii}}\right]_{\hat{u}_{min,i}}^{\hat{u}_{max,i}}$;
            \item $\Delta \leftarrow \theta - \hat{u}_{t,i}$;
            \item $\sigma \leftarrow \sigma + \Delta^2$;
            \item $\hat{u}_{t,i} \leftarrow \theta$;
            \item $V \leftarrow V + \Delta \hat{B}_{\cdot,i}$;
        \end{enumerate}
    \item \textbf{end}.
    \end{enumerate}
    \vspace*{.1cm}\hrule\vspace*{.1cm}
    \textbf{Output}: $\hat u_{t}, \sigma$. 
\end{algorithm}

\subsection{Preconditioning}
Preconditioning is a common heuristic for improving the computational performance of first-order methods. The optimal design of preconditioners has been studied for several decades, but such a computation is often more complex than the original problem and may become prohibitive if it must be executed at runtime. Diagonal
scaling is a heuristic preconditioning that is very simple and often beneficial~\cite{giselsson2014diagonal,takapoui2016preconditioning}. In this paper, we propose
to make the change of state variables $\bar{x} =E\hat{x}$, where $E \in \mathcal{R}^{\hat{n}_x \times \hat{n}_x}$ is a diagonal matrix whose
$i$th entry is
\begin{equation}
E_{i,i}=\sqrt{Q_{i,i}+\hat{A}_{\cdot,i}^{\prime}\hat{A}_{\cdot,i}}
\end{equation}
and replace the prediction model $\hat{x}_{t+1}=\hat{A}\hat{x}_t+\hat{B}\hat{u}_t$ by
\[
\bar{x}_{t+1}=\bar{A}\bar{x}_t+\bar{B}\hat{u}_t
\]
where $\bar{A} = E\hat{A}E^{-1}$ and $\bar{B}=E\hat{B}$. The weight matrix $Q$ and constraints $[\hat{x}_{\rm min},\hat{x}_{\rm max}]$ are scaled accordingly by setting $\bar{Q} = E^{-1}QE^{-1}$ and $\bar{x}_{\rm min}=E^{-1}\hat{x}_{\rm min}$, $\bar{x}_{\rm max}=E^{-1}\hat{x}_{\rm max}$.

\subsection{Efficient coupling scheme between CD and AL method}
We are now ready to couple CD and AL to solve the posed MPC problem~\eqref{MPC}
efficiently. We first note that updating $u_t$ and $x_{t+1}$ for all $t$ involves computing a similar temporary vector $V$ in Procedure~\ref{DetailedRCCD}. As $V$ is in fact the next update of the dual vector $\Lambda$
in Algorithm~\ref{AALM}, we modify Procedure~\ref{DetailedRCCD} as shown in Procedure~\ref{CouplingRCCD}. The overall solution method described in the previous subsections  is summarized in Algorithm~\ref{ccdALM}, that we call CDAL. Note that the main update of the Lagrangian variables in Algorithm~\ref{ccdALM} is placed early in Step~\ref{algo:ccdALM:dual}, unlike in Algorithm~\ref{AALM},due to the use of the proposed efficient coupling scheme. The AL (outer) iterations are executed for maximum $N_{\rm out}$ iterations, the CD (inner) iterations for at most $N_{\rm in}$ iterations. The tolerances $\epsilon_{\rm out}$ and $\epsilon_{\rm in}$ are used to stop the outer and inner iterations, respectively. Algorithm~\ref{ccdALM} is matrix-free and library-free, and we could
implement it in 90 lines of C code.

\floatname{algorithm}{Procedure}
\begin{algorithm}[t]
    \caption{Modified Procedure~\ref{DetailedRCCD} to efficiently couple CD and AL}\label{CouplingRCCD}
    \textbf{Input}: $\lambda_t$, $\hat{u}_{t}$;
    MPC settings $\hat{A}$, $\hat{B}$, $R$, $\hat u_{\min}$, $\hat u_{\max}$; parameter $\rho>0$; update amount $\sigma\geq 0$.
    \vspace*{.1cm}\hrule\vspace*{.1cm}
    \begin{enumerate}[label*=\arabic*., ref=\theenumi{}]
        \item \textbf{for} $i=n_u,\ldots, 1$ \textbf{do}
        \begin{enumerate}[label=\theenumi{}.\arabic*., ref=\theenumi{}.\arabic*]
            \item $s \leftarrow  \frac{1}{\rho}R_{i,\cdot} \hat{u}_{t} + (\hat{B}_{\cdot, i})^{\prime} \lambda_t$;
            \item $\theta \leftarrow  \left[\hat{u}_{t,i}-\frac{s}{\frac{1}{\rho}R_{ii}+(\hat{B}^{\prime}\hat{B})_{ii}}\right]_{\hat{u}_{min,i}}^{\hat{u}_{max,i}}$;
            \item $\Delta \leftarrow \theta - \hat{u}_{t,i}$;
            \item $\sigma \leftarrow \sigma + \Delta^2$;
            \item $\hat{u}_{t,i} \leftarrow \theta$;
            \item $\lambda_{t} \leftarrow \lambda_t + \Delta \cdot \hat{B}_{\cdot,i}$;
        \end{enumerate}
    \item \textbf{end}.
    \end{enumerate}
    \vspace*{.1cm}\hrule\vspace*{.1cm}
    \textbf{Output}: $\hat u_{t}, \lambda_t, \sigma$.
\end{algorithm}

\floatname{algorithm}{Algorithm}
\begin{algorithm}
    \caption{Accelerated reverse cyclic CDAL algorithm for linear (or linearized) MPC}\label{ccdALM}
    \textbf{Input}: primal/dual warm-start $U = \{\hat{u}_{0}, \hat{u}_{1}, \cdots, \hat{u}_{T-1}\}$, $X =\{\hat{x}_{0}, \hat{x}_{1}, \cdots, \hat{x}_{T}\}$, $\Lambda^{-1}$ $=$ $\Lambda^0$ $=$ $\{\lambda_{0}$, $\lambda_{1}$, $\cdots$, $\lambda_{T-1}\}$; MPC settings $\{\hat{A},$ $\hat{B},$ $\hat{C},$ $W_y$, $W_u$ ,$W_{\Delta u}$, $\Delta u_{\min}$, $\Delta u_{\max}$, $u_{\min}$, $u_{\max}$, $x_{\min}$, $x_{\max}\}$; Algorithm settings $\{\rho, N_{\rm out},N_{\rm in}\,\epsilon_{\rm out}, \epsilon_{\rm in}\}$
    \vspace*{.1cm}\hrule\vspace*{.1cm}
    \begin{enumerate}[label*=\arabic*., ref=\theenumi{}]
        \item Obtain preconditioned $\bar{X}=\{\bar{x}_{0},\cdots, \bar{x}_{T}\}$, $\bar{A}$, $\bar{B}$, $\bar{Q}$, $\bar{x}_{\rm min}, \bar{x}_{\rm max}$ according to Section III.C
        \item $\alpha_1\leftarrow1$; $\hat{\Lambda}^{0}\leftarrow\Lambda^{0}$;
        \item \textbf{for} $k=1,2,\cdots, N_{\rm out}$ \textbf{do}
        \begin{enumerate}[label=\theenumi{}.\arabic*., ref=\theenumi{}.\arabic*]
            \item\textbf{for} $t=0,\ldots,T-1$ \textbf{do}\label{algo:ccdALM:dual}
            \begin{enumerate}[label=\theenumii{}.\arabic*., ref=\theenumii{}.\arabic*]
                \item $\lambda^{k}_{t} = \hat{\lambda}^{k-1}_{t} + \bar{A} \bar{x}_{t} + \bar{B} \hat{u}_{t} - \bar{x}_{t+1}$;
            \end{enumerate}
            \item \textbf{for}  $k_{in}=1,2,\cdots, N_{\rm in}$  \textbf{do}         
            \begin{enumerate}[label=\theenumii{}.\arabic*., ref=\theenumii{}.\arabic*]
                \item $U, \bar{X}, \sigma \leftarrow$  Procedure~\ref{RCCD} with use of Procedure 4;
                \item \textbf{if} $\sigma \leq \epsilon_{\rm in}$ \textbf{break} the loop;
            \end{enumerate}            
            \item \textbf{if} $\|\Lambda^{k}-\hat{\Lambda}^{k-1}\|_2^2\leq \epsilon_{\rm out}$ \textbf{stop};
            \item $\alpha_{k+1} \leftarrow \frac{1+\sqrt{1+4\alpha_k^2}}{2}$;
            \item $\hat{\Lambda}^{k} \leftarrow \Lambda^k + \frac{\alpha_k-1}{\alpha_{k+1}} (\Lambda^k-\Lambda^{k-1})$;
        \end{enumerate}
    \item Recover $X$ from $\bar{X}$
    \item \textbf{end}.
    \end{enumerate}
    \vspace*{.1cm}\hrule\vspace*{.1cm}
    ~~\textbf{Output}: $U, X, \Lambda$
\end{algorithm}

\section{Numerical Examples}
\label{sec:examples}
We test the performance of the CDAL solver against other solvers in two numerical experiments. The first one is the ill-conditioned AFTI-16 control problem \cite{KAS88, bemporad1997nonlinear} based on LTI-MPC, used in the Model Predictive Control Toolbox for MATLAB~\cite{bemporad2004model}. The main goals of this experiment include investigating whether our proposed simple heuristic preconditioner, reverse cyclic rule, and Nesterov's acceleration scheme are helpful, and provide a detailed comparison with other solvers. The second experiment demonstrates the benefits of the construction-free property in LPV-MPC of a CSTR~\cite{seborg2010process}, in which
the prediction model is obtained by linearizing a nonlinear model of the process at each sample step. The reported simulation results were obtained on a MacBook Pro with 2.7~GHz 4-core Intel Core i7 and 16GB RAM. Algorithm~\ref{ccdALM} is executed in MATLAB via a C-mex interface.

\subsection{AFTI-16 Benchmark Example}
The open-loop unstable linearized AFTI-16 aircraft model reported in \cite{KAS88,bemporad1997nonlinear} 
is 
\[
\left\{\begin{aligned}
\dot{x} =&{\footnotesize\left[\begin{array}{cccc}
-0.0151 & -60.5651 & 0 & -32.174 \\
-0.0001 & -1.3411 & 0.9929 & 0 \\
0.00018 & 43.2541 & -0.86939 & 0 \\
0 & 0 & 1 & 0
\end{array}\right]} x\\&+{\footnotesize\left[\begin{array}{cc}
-2.516 & -13.136 \\
-0.1689 & -0.2514 \\
-17.251 & -1.5766 \\
0 & 0
\end{array}\right] }u \\
y =&{\footnotesize\left[\begin{array}{llll}
0 & 1 & 0 & 0 \\
0 & 0 & 0 & 1
\end{array}\right]x}
\end{aligned}\right.
\]
The model is sampled using zero-order hold every 0.05~s. The input constraints are $|u_i| \leq 25^{\circ},i = 1, 2$, the output constraints are $-0.5\leq y_1 \leq 0.5 $ and $-100 \leq y_2 \leq 100$. The control goal is to make the pitch angle $y_2$ track a reference signal $r_2$. In designing the MPC controller we
take $W_y = \diag$([10,10]), $W_u = 0$, $W_{\Delta u}= \diag$([0.1, 0.1]), and the
prediction horizon is $T=5$.  

To investigate the effects of the three techniques (reverse cyclic rule, acceleration, and preconditioning) that we have introduced to improve the efficiency of the CDAL algorithm, we performed closed-loop simulations on eight schemes with fixed $\rho=1$. These are:
0-CDAL, the basic scheme, without acceleration and reverse cyclic rule; {R-CDAL}, the scheme with the Reverse cyclic rule; {A-CDAL}, the Accelerated scheme; {AR-CDAL}, the Accelerated scheme with the Reverse cyclic rule, and their respective schemes with preconditioner, namely {P-0-CDAL}, {P-R-CDAL}, {P-A-CDAL}, and finally {CDAL}, that includes all 
the proposed techniques. The stopping criteria are 
defined by $\epsilon_{\rm in}=10^{-6}$, $\epsilon_{\rm out}=10^{-4}$, and 
$N_{\rm out}$, $N_{\rm in}$ are set to the large enough value 5000 in order
to guarantee good-quality solutions. 
 
\begin{figure}
        \hspace*{-1em}\includegraphics[width=1.1\columnwidth]{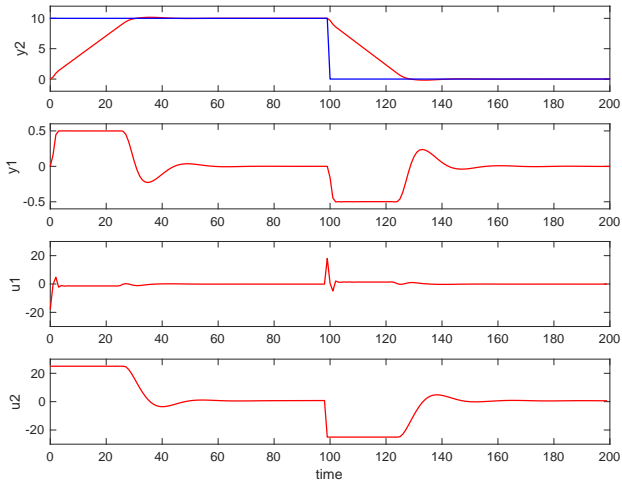} 
        \caption{Linear AFTI-16 closed-loop performance}
        \label{Fig1} 
\end{figure}

The computational load associated with the above schemes is listed in Table~\ref{tab1}, in which the last column represents the closed-loop performance, which is the average value $\frac{1}{T}\sum_{t=0}^{T-1} \left\|W_y \left(y_{t+1}-r_{t+1}\right)\right\|_{2}^{2} + \left\|W_u  \left(u_{t+1}-u^r_{t+1}\right)\right\|_{2}^{2} + \left\|W_{\Delta u}  \Delta u_{t}\right\|_{2}^{2}$ of the MPC cost 
over the duration $T$ of the closed-loop simulation and is almost the same for all schemes. 
The associated closed-loop trajectories are reported in Figure~\ref{Fig1}, which shows that the pitch angle correctly tracks the reference signal from $0^{\circ}$ to $10^{\circ}$ and then back to $0^{\circ}$, and that both the input and output constraints are satisfied.

Since each MPC execution requires different numbers of inner and outer iterations, the average (``avg'') and maximum (``max'') number of iterations (or CPU time) are computed over the entire closed-loop execution. It can be observed that the maximum and average number of inner-loop iterations of R-CDAL are smaller than that of {0-CDAL} (especially the maximum number), while their outer-loop iterations are almost the same, which shows that the reverse cyclic rule provides a significant improvement. Although {A-CDAL} has fewer outer-loop iterations, it has more inner-loop iterations than {0-CDAL} on average. It therefore does not result in a significant reduction in
total computation time. We can see that {AR-CDAL} achieves fewer iterations both in the inner loop and outer loop and has better average and worst-case computation performance. It can also be seen from Table~\ref{tab1} that preconditioning significantly reduces the number of outer-loop iterations.

\begin{table}[!htbp]
\caption{Computational performance of different schemes}
\centering
\begin{tabular}{lrrrrrrr}
\toprule
method  & \multicolumn{2}{c}{sum of inner iters}& \multicolumn{2}{c}{outer iters}& \multicolumn{2}{c}{time (ms)} & cost\\
&avg&max&avg&max&avg&max&
\\\midrule
{0-CDAL} & 8577&79615 & 339&2104 & 4.9&55.3 &42.3\\
{R-CDAL}& 7298&72693 & 340&2103 & 4.3&53.2 &42.5\\
{A-CDAL} & 7437&57026 & 45&297 & 4.0&41.1 &42.5\\
{AR-CDAL} & 6207&51884 & 44&205 & 3.8&39.5 &42.5\\
\bottomrule
{P-0-CDAL} & 3467&13386 & 33&171 & 2.1&11.4 &42.5\\
{P-R-CDAL}& 1757&13430 & 33&171 & 1.0&10.9 &42.5\\
{P-A-CDAL} & 3299&12161 & 13&60 & 1.7&9.7 &42.5\\
\textbf{CDAL} & 1543&12508 & 13&60 & 0.85&9.5&42.5\\
\bottomrule
\end{tabular}
\label{tab1}
\end{table}

Next, we investigate the effect on computation efficiency of parameter $\rho$,
that we expect to tend to trade off feasibility versus optimality.
In particular, we expect larger values of $\rho$ to favor feasibility, i.e., provide more inner-loop
iterations and less outer-loop iterations, and vice versa. The computational performance results obtained by performing closed-loop simulations using the final \textbf{CDAL} algorithm for different values of $\rho$ between 0.01 and 1 are listed in Table~\ref{tab2}. When the parameter value is between 0.01 and 0.1, the {CDAL} algorithm has very similar computational burden.

To further illustrate the efficiency of CDAL, Table~\ref{tab2} also lists the results obtained by using other solvers. Here the fastMPC solver is also a construction-free solver which provides a free C-mex code. We also made comparison with the $\mu$AO-MPC solver v1.0.0-beta \cite{zometa2013muao}, which is based on an augmented Lagrangian method together with Nesterov’s gradient method. The $\mu$AO-MPC differs from CDAL in the way the subproblems are solved, and the outer loop not involving an acceleration scheme. The state-of-the-art first-order method for QP, the OSQP solver v0.6.2~\cite{SBGBB20}, and MATLAB's built-in QP solver (quadprog) are also used for comparison. For a fair comparison, each solver setting is chosen to at least ensure each shares the same objective cost and constraint violation. When the parameter $\rho$ of the {CDAL} is 0.01, the CDAL is faster than the other solvers. Regarding the $\mu$AO-MPC, OSQP and quadprog solver, we split between QP problem construction time (including the required matrix factorizations) and pure solution time. Note that in this case, the controller is LTI-MPC, and hence the MPC problem construction and matrix factorizations required by these non-construction-free solvers can be performed offline. On the other hand, in case of LPV-MPC problems the total computation time would be spent online and the embedded code would also include routines for problem construction and matrix factorization functions. Instead, CDAL does not require any construction nor factorizations, thus making the solver very lean and fast also in a time-varying
MPC setting, as investigated next.

\begin{table}[!htbp]
\caption{Computational load of CDAL with different values of $\rho$ and comparison with other solvers}
\centering
\begin{threeparttable}
\begin{tabular}{ccccc}
\toprule
Solver  & solver setting& \multicolumn{2}{c}{time (ms)} & cost\\
& & avg&max&
\\\midrule
CDAL&$\rho=1  $ & 0.85&9.5 & 42.561\\
    &$\rho=0.5$ & 0.72&7.1 & 42.590\\
    &$\rho=0.2$ & 0.53&4.2 & 42.612\\
    &$\rho=0.1$ & 0.47&3.8 & 42.619\\
    &$\rho=0.05$& 0.42&3.3 & 42.618\\
    &$\rho=0.01$& 0.41&3.2 & 42.618\\
\bottomrule
FastMPC & $maxit=5,k= 0.1$ & 0.54&4.2 & 42.627\\
\bottomrule
$\mu$AO-MPC    & $\mu=0.05$ &7.0\tnote{*}&68.1\tnote{*} & 42.627\\
        &   in_iter=100,ex_iter=100        &8\tnote{**} & 69\tnote{**}\\
\bottomrule
OSQP    & $N=5000,\epsilon=10^{-6}$& 0.6\tnote{*}&10.1\tnote{*} & 42.627 \\
        &           &1.5\tnote{**} & 13.8\tnote{**}\\
\bottomrule
quadprog& default & 10.3\tnote{*}&20.6\tnote{*} & 42.622\\
        &      &11\tnote{**} & 22\tnote{**}\\

\bottomrule
\end{tabular}
\begin{tablenotes}
    \footnotesize
    \item[*]: pure solution time, without including matrix factorization
    \item[**]: total time (MPC construction + solution)
\end{tablenotes}
\end{threeparttable}
\label{tab2}
\end{table}

\subsection{Nonlinear CSTR Example}
To illustrate the performance of CDAL when the linear MPC formulation~\eqref{MPC} changes at runtime we consider the 
control of the CSTR system \cite{seborg2010process}, described by the continuous-time nonlinear
model
\begin{equation}
\begin{array}{rcl}
\frac{dC_A}{dt} &=& C_{A,i}-C_A-k_0 e^{\frac{-EaR}{T}}C_A \\
\frac{dT}{dt} &=& T_{i} + 0.3 T_c - 1.3T + 11.92 k_0 e^{\frac{-EaR}{T}}C_A\\
y&=&C_A
\end{array}
\label{eq:CSTR}
\end{equation}
where $C_A$ is the concentration of reagent A, $T$ is the temperature of the reactor, $C_{A,i}$ is the inlet feed stream concentration, which is assumed to have the constant value $10.0$ kgmol/m$^3$. The disturbance comes from the inlet feed stream temperature $T_{i}$, which has fluctuations represented by $T_{i} = 298.15 + 5 \sin(0.05 t)$ $K$. The manipulated variable is the coolant temperature $T_c$. The constants $k_0 = 34930800$ and $EaR = -5963.6$ (in MKS units). The reactor's initial state is at a low conversion rate, with $C_A = 8.57$ kgmol/m$^3$, $T = 311$ K. The goal is to adjust the reactor state to a high reaction rate with $C_A = 2$ kgmol/m$^3$, which
is a quite large condition. The controller manipulates the coolant temperature $T_c$ to track a concentration reference as well as reject the measured disturbance $T_{i}$. Due to its nonlinearity, the
model in~\eqref{eq:CSTR} is linearized online at each sampling step:
\[
\frac{dx}{dt} \approx f(x_{t},u_{t-1},p) + \left.\frac{\partial f}{\partial x}\right|_{x_t, u_{t-1}, p}\hspace*{-2em}(x-x_{t}) + \left.\frac{\partial f}{\partial u}\right|_{x_t, u_{t-1}, p}
\hspace*{-2em}(u-u_{t-1})
\]
where $f(x,u,p)$ is the mapping defined in~\eqref{eq:CSTR} for $x=[C_A\ T]'$, $u=T_c$, $p=[C_{A,i}\ T_i]'$.
By setting $A_c=\left.\frac{\partial f}{\partial x}\right|_{x_t, u_{t-1}, p}, B_c=\left.\frac{\partial f}{\partial u}\right|_{x_t, u_{t-1}, p}, e_c = f(x_{t},u_{t-1},p) - A_t x_t - B_t u_{t-1}$, we get the following linearized continuous-time model
\[
\frac{d}{dt}x = A_c x + B_c u + e_c
\]
We use the forward Euler method with sampling time $T_s = 0.5$ minutes to obtain the following discrete-time model
\[
x_{t+1} = A_d x_t + B_d u_t+ e_d
\]
where $A_d = I + T_s A_c, B_d = T_s B_c, e_d = T_s e_c$. Although held constant over the prediction horizon, clearly matrices $A_d, B_d$ and the offset term $e_d$ change at runtime, which makes the controller an LPV-MPC. Regarding the performance index, we choose weights $W_y =1$, $W_u = 0$, $W_{\Delta u} = 0.1$. The physical limitation of the coolant jacket is that its rate of change $\Delta T_c$ is subject to the constraint $[-1,1]$~K when considering the sampling time $T_s=0.5$ minutes. The prediction horizon is $T=10$ steps. 

\begin{figure}
        \hspace*{-1em}\includegraphics[width=1.1\columnwidth]{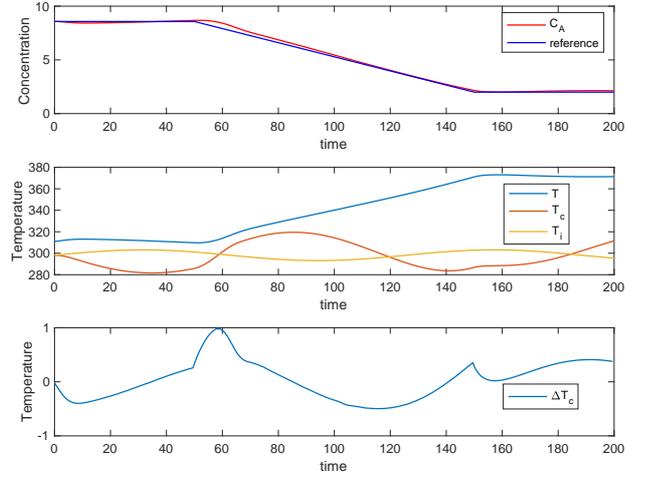} 
        \caption{Nonlinear CSTR closed-loop performance}
        \label{Fig2} 
\end{figure}

We compare again CDAL with fastMPC, FGAL, OSQP, and quadprog solvers in the LPV-MPC setting described above. CDAL is run with $\epsilon_{\rm in}=10^{-6}$, $\epsilon_{\rm out}=10^{-4}$, $\rho=0.01$, and $N_{\rm out}=N_{\rm in}=5000$. For a fair comparison, each solver setting is chosen to at least ensure each shares the same objective cost and constraint violation. The closed-loop simulation results of CDAL and other solvers almost coincide and are plotted in Figure~\ref{Fig2}, from which it can be seen that $C_A$ tracks the reference signal well, and the fluctuation of $T_i$ is effectively suppressed. The computational load and closed-loop performance associated with CDAL and other solvers are reported in Table~\ref{tab3}. In this successive linearization-based MPC example, we found that the problem-construction time has a comparable computation time to the problem-solving time from the results of non-construction-free solvers. If we only compare the solution time, CDAL is faster than other solvers except for OSQP, but in fact the MPC construction time must be included for comparison, which leads to CDAL being faster than OSQP. Because of the construction-free, matrix-free, and library-free features, CDAL has an advantage in industrial embedded deployment when the optimization problem associated with MPC is constructed online and this operation has a cost that is comparable to the solution time.
\begin{table}[!htbp]
\caption{Computational performance of CDAL and other solvers}
\centering
\begin{threeparttable}
\begin{tabular}{ccccc}
\toprule
Solver  & solver setting& \multicolumn{2}{c}{time (ms)} & cost\\
& & avg&max&
\\\midrule
CDAL  &$\rho,\epsilon_{in},\epsilon_{out}=0.01,10^{-6},10^{-4}$ & 0.3 &0.6 & 0.02202\\
\hline
FastMPC &maxit=5,$k=0.1$  & 0.5 & 7.2&0.030170 \\
\hline
$\mu$AO-MPC &$\mu=0.01$  & 1.4\tnote{*} & 10.1\tnote{*} &0.02202\\
&in_iter=100,ex_iter=10 & 2.1\tnote{**} & 15.2\tnote{**}\\
\hline
OSQP  & default & 0.15\tnote{*} & 0.37\tnote{*} & 0.02219\\
      & & 0.6\tnote{**} & 5.5\tnote{**}\\
\hline
quadprog & default & 1.6\tnote{*} & 9.7\tnote{*} & 0.02219 \\
& & 1.8\tnote{**} & 13.3\tnote{**}\\
\bottomrule
\end{tabular}
\begin{tablenotes}
    \footnotesize
    \item[*]: solution time
    \item[**]: MPC construction time + solution time
\end{tablenotes}
\end{threeparttable}
\label{tab3}
\end{table} 

\section{Conclusion}
\label{sec:conclusions}
This paper has proposed a construction-free, matrix-free, and library-free MPC solver, based on a cyclic coordinate-descent method in the augmented Lagrangian framework. We showed that the method is efficient and competes with other existing methods, thanks to the use of a reverse cyclic rule, Nesterov’s acceleration, a simple heuristic preconditioner, and an efficient coupling scheme. Compared to many QP solution methods proposed in the literature, CDAL avoids constructing the QP problem, which makes it particularly appealing for some scenarios in which its online construction is required and has a comparable computation time to solving itself. 

The proposed algorithm can be immediately extended to handle linear time-varying systems, in which
the plant-model and/or cost-function matrices are allowed to vary over the prediction horizon. 
Future research will investigate the use of CDAL to solve nonlinear MPC problems and data-driven MPC formulations in which the model is adapted online by recursive system identification.

\bibliographystyle{unsrt}
\bibliography{refs} 
\end{document}